\documentclass{article}

\usepackage[margin=3cm]{geometry}

\usepackage{latexsym,amssymb,amsmath,amsthm,amsfonts,amscd,mathtools}
\usepackage{enumerate}
\usepackage{graphicx}
\usepackage{tikz}
\usepackage{tikz-3dplot}
\usetikzlibrary{3d,calc}
\usetikzlibrary{matrix,arrows}
\usepackage{multicol}
\usepackage{multirow}
\usepackage{tabularx}
\usepackage{listings}
\usepackage{color}
\usepackage{hyperref}
\usepackage{kbordermatrix}
\newcommand{\vd}{\vdots}
\newcommand{\cd}{\cdots}
\newcommand{\dd}{\ddots}

\usepackage{color}
\definecolor{deepblue}{rgb}{0,0,0.5}
\definecolor{deepred}{rgb}{0.6,0,0}
\definecolor{deepgreen}{rgb}{0,0.5,0}
\definecolor{cof}{RGB}{219,144,71}
\definecolor{pur}{RGB}{186,146,162}
\definecolor{greeo}{RGB}{91,173,69}
\definecolor{greet}{RGB}{52,111,72}

\newtheorem{theorem}{Theorem}
\newtheorem{lemma}[theorem]{Lemma}

\newtheorem{corollary}[theorem]{Corollary}

\theoremstyle{definition}

\newtheorem{example}[theorem]{Example}

\newcommand{\diag}{\operatorname{diag}}
\newcommand{\SNF}{\operatorname{SNF}}
\DeclareMathOperator\dist{dist}
\def \rank {\operatorname{rank}}

\title{The Smith normal form of distance matrices of high dimensional trees}

\author{Carlos A. Alfaro\thanks{
{\tt alfaromontufar@gmail.com}, 
Banco de M\'exico, 
Ciudad de M\'exico, Mexico} 
\and Jesús Uriel Medrano\thanks{
{\tt jesusurielmedrano@gmail.com}, 
Departamento de Matemáticas, Facultad de Ciencias, 
Universidad Nacional Autónoma de México, 
Ciudad de M\'exico, México.} 
\and Iván Téllez Téllez\thanks{
{\tt ivan.tellez@uaslp.mx}, 
Facultad de Economía,
Universidad Autónoma de San Luis Potosí, 
San Luis Potosí, México.}}

\begin{document}

\maketitle


\begin{abstract}
Graham-Lovász-Pollak \cite{GL,GP} obtained the celebrated formula $$\det({\sf D}(T_{n+1}))=(-1)^nn2^{n-1},$$ for the determinant of the distance matrix ${\sf D}(T_{n+1})$ for any tree $T_{n+1}$ with $n+1$ vertices.
Later, Hou and Woo \cite{HW} extended this formula to the Smith normal form (SNF) obtaining that $\SNF({\sf D}(T_{n+1}))={\sf I}_2\oplus 2{\sf I}_{n-2}\oplus [2n]$, for any tree $T_{n+1}$ with $n+1$ vertices.

A $k$-{\it tree} is either a complete graph on $k$ vertices or a graph obtained from a smaller $k$-tree by adjoining a new vertex together with $k$ edges connecting it to a $k$-clique.
If $\tau$ and $\tau'$ are $d$-cliques in a $k$-tree $T$, a $d$-{\it walk} between $\tau$ and $\tau'$ is a finite sequence $\tau_1\sigma_1\tau_2\sigma_2\cdots\tau_l$, where $\tau_1=\tau$, $\tau_l=\tau'$, and the $d$-cliques $\tau_i$ and $\tau_{i+1}$ are incident to the same $(d+1)$-clique $\sigma_i$.
For $d\in\{1,\dots,k\}$, the $d$-{\it distance} from the $d$-cliques $\tau$ and $\tau'$ is the number of $(d+1)$-cliques in a minimum $d$-walk from $\tau$ and $\tau'$, and is denoted by $\dist^d(\tau,\tau')$.
Let $c_d$ denote the number of $d$-cliques in the $k$-tree $T$. 
Then the $d$-distance matrix ${\sf D}^d(T)$ of the $k$-tree $T$ is the $c_d\times c_d$ matrix, indexed by the $d$-cliques of $T$, such that the $(i,j)$-entry is $0$ if $i=j$, and $\dist^d(\tau_i,\tau_j)$ otherwise.
Here, we show that, for $k$ and $n$ fixed, the SNF of the $k$-distance matrix is the same for any $k$-tree with $n$ vertices.
Specifically, for any $k$-tree $T_{n}$ with $n$ vertices such that $n\geq k+2$, the Smith normal form of ${\sf D}^{k}(T_{n})$ is 
$${\sf I}_{(k-1)(n-k)+2}\oplus (k+1){\sf I}_{n-k-2}\oplus [k(k+1)(n-k)],$$
which extends Graham-Lovász-Pollak and Hou-Woo results.
\end{abstract}



\section{Introduction}
\begin{figure}[!ht]
    \centering
    \begin{tabular}{c@{\extracolsep{1cm}}c}
        \begin{tikzpicture}[scale=3,thick]
            \tikzstyle{every node}=[draw,minimum width=0pt, inner sep=2pt, circle, fill=white]
            \coordinate (A1) at (0,0);
            \coordinate (A2) at (0.6,0.2);
            \coordinate (A3) at (1,0);
            \coordinate (A4) at (0.4,-0.2);
            \coordinate (B1) at (0.5,0.5);
            \coordinate (B2) at (0.5,-0.5);
            
            \draw (B1) -- (A4) -- (A1) -- (A2) -- (A4) -- (A3) -- (B2) -- (A4) -- (B1) -- (A1) -- (B2); 
    
            \draw (A1) node[label=left:{\footnotesize $v_0$}] {};
            \draw (A2) node[label=right:{\footnotesize $v_1$}] {};
            \draw (A3) node[label=0:{\footnotesize $v_2$}] {};
            \draw (A4) node[label=100:{\footnotesize $v_3$}] {};
            \draw (B1) node[label=north:{\footnotesize $v_4$}] {};
            \draw (B2) node[label=south:{\footnotesize $v_5$}] {};
        \end{tikzpicture}
         & 
        \begin{tikzpicture}[thick,scale=3]
            \coordinate (A1) at (0,0);
            \coordinate (A2) at (0.6,0.2);
            \coordinate (A3) at (1,0);
            \coordinate (A4) at (0.4,-0.2);
            \coordinate (B1) at (0.5,0.5);
            \coordinate (B2) at (0.5,-0.5);

            \draw[fill=cof,opacity=0.6] (A1) -- (A4) -- (B1);
            \draw[fill=pur,opacity=0.6] (A1) -- (A4) -- (B2);
            \draw[fill=greeo,opacity=0.6] (A1) -- (A2) -- (A4);
            \draw[fill=blue!50,opacity=0.6] (A3) -- (A4) -- (B2);
            \draw (B1) -- (A1) -- (B2) -- (A3);
        
            \draw (A1) node[label=left:{\footnotesize $v_0$}] {};
            \draw (A2) node[label=right:{\footnotesize $v_1$}] {};
            \draw (A3) node[label=0:{\footnotesize $v_2$}] {};
            \draw (A4) node[label=100:{\footnotesize $v_3$}] {};
            \draw (B1) node[label=north:{\footnotesize $v_4$}] {};
            \draw (B2) node[label=south:{\footnotesize $v_5$}] {};
        \end{tikzpicture}
    \end{tabular}
    
    \caption{A 2-tree (left) and a simplicial complex realization (right).}
    \label{fig:realization of a 2 tree}
\end{figure}
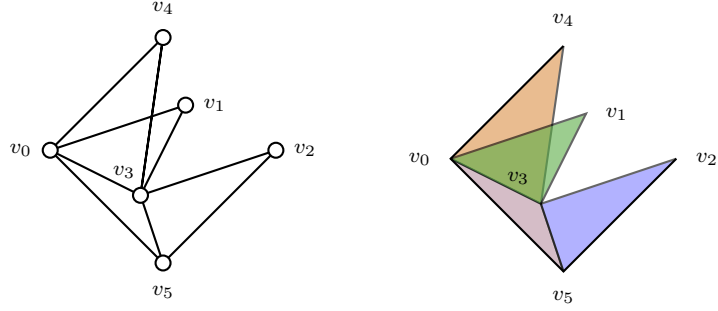

A $k$-{\it clique} is a complete subgraph on $k$ vertices.
The concept of $k$-trees may be defined recursively: a $k$-{\it tree} is either a complete graph on $k$ vertices or a graph obtained from a smaller $k$-tree by adjoining a new vertex together with $k$ edges connecting it to a $k$-clique.
These concepts can be analogously defined in terms of simplicial complexes and simplices, for instance in Figure~\ref{fig:realization of a 2 tree} we observe a simplicial realization of a 2-tree.
For further reading in simplicial complexes, we refer the reader to the book of Munkres \cite{MR755006}.

Harary and Palmer \cite{MR228355,MR357214} counted unlabeled 2-trees in 1968.
The enumeration of unlabeled $k$-trees for $k > 2$ was a long-standing unsolved problem
until the recent solution by Gainer-Dewar \cite{MR3007180}, using the theory of combinatorial
species.
In \cite{MR3213312}, unlabeled $k$-trees are enumerated by first coloring the vertices using $k+1$ colors in such a way that the vertices of each $(k+1)$-clique use all $k+1$ colors.
See the sequences \href{https://oeis.org/A054581}{A054581} and \href{https://oeis.org/A370770}{A370770} in The On-Line Encyclopedia of Integer Sequences \cite{oeis}.
In Appendix \ref{sec:generation of k trees}, we include {\tt sagemath} code that generates the $k$-trees with up to $n$ vertices.

Let $T$ be a $k$-tree.
For $d\in \{1,\dots,k\}$, two $d$-cliques $\tau$ and $\tau'$ in $T$ are {\it adjacent} if they belong to the same $(d+1)$-clique $\sigma$, in such situation $\tau$ and $\tau'$ are {\it incident} to $\sigma$. 
In this way, if $\tau$ and $\tau'$ are $d$-cliques, a $d$-{\it walk} between $\tau$ and $\tau'$ is a finite sequence $\tau_1\sigma_1\tau_2\sigma_2\cdots\tau_l$, where $\tau=\tau_1$, $\tau'=\tau_l$, and $\tau_i$ and $\tau_{i+1}$ are incident to the same $(d+1)$-clique $\sigma_i$.
Therefore, the $d$-{\it distance} from the $d$-cliques $\tau$ and $\tau'$ is the number of $(d+1)$-cliques in a minimum $d$-walk from $\tau$ and $\tau'$, and is denoted by $\dist^d(\tau,\tau')$.
Note that there exists a $d$-walk between any pair of $d$-cliques in any $k$-tree $T$.
Let $c_d$ denote the number of $d$-cliques in $T$. 
Then the $d$-{\it distance matrix} ${\sf D}^d(T)$ of the $k$-tree $T$ is the $c_d\times c_d$ matrix, indexed by the $d$-cliques of $T$, such that
\[
{\sf D}^d(T)_{i,j}=
\begin{cases}
    0 & \text{if } i=j\\
    \dist^d(\sigma_i,\sigma_j) & \text{otherwise.}
\end{cases} 
\]
For instance, if $T$ is the 2-tree shown in Figure~\ref{fig:realization of a 2 tree}, then
\[
{\sf D}^1(T)=
\begin{bmatrix}
0 & 1 & 2 & 1 & 1 & 1 \\
1 & 0 & 2 & 1 & 2 & 2 \\
2 & 2 & 0 & 1 & 2 & 1 \\
1 & 1 & 1 & 0 & 1 & 1 \\
1 & 2 & 2 & 1 & 0 & 2 \\
1 & 2 & 1 & 1 & 2 & 0 \\
\end{bmatrix}
\]
and
\[
{\sf D}^2(T) =
	    \kbordermatrix{
	           &  01 &  03 &  04 &  05 &  13 &  23 &  25 &  34 &  35 \\
            01 &   0 &   1 &   2 &   2 &   1 &   3 &   3 &   2 &   2 \\
            03 &   1 &   0 &   1 &   1 &   1 &   2 &   2 &   1 &   1 \\
            04 &   2 &   1 &   0 &   2 &   2 &   3 &   3 &   1 &   2 \\
            05 &   2 &   1 &   2 &   0 &   2 &   2 &   2 &   2 &   1 \\
            13 &   1 &   1 &   2 &   2 &   0 &   3 &   3 &   2 &   2 \\
            23 &   3 &   2 &   3 &   2 &   3 &   0 &   1 &   3 &   1 \\
            25 &   3 &   2 &   3 &   2 &   3 &   1 &   0 &   3 &   1 \\
            34 &   2 &   1 &   1 &   2 &   2 &   3 &   3 &   0 &   2 \\
            35 &   2 &   1 &   2 &   1 &   2 &   1 &   1 &   2 &   0 \\
         }.
\]
In Appendix \ref{sec:k distance matrix}, we give the {\tt sagemath} function {\tt d\_distance\_matrix(G,d)} that computes the $d$-distance matrix of a $k$-tree $G$. 
In particular, for the $2$-tree shown in Figure~\ref{fig:realization of a 2 tree} it can be called with the following code:
\begin{lstlisting}
S = SimplicialComplex([[0,1,3],[0,3,4],[0,3,5],[2,3,5]])
G = S.graph()
G.show()
d_distance_matrix(G,1)
d_distance_matrix(G,2)
\end{lstlisting}

Note that ${\sf D}^1(T)$ coincides with the usual distance matrix of the underlying graph of the $k$-tree $T$.
For more information on the distance matrix of a graph, we refer the reader to \cite{MR3231823,MR4482731,MR4239911}.
In the following, ${\sf I}_n$ and ${\sf J}_n$ will denote the identity and the all-ones matrices of order $n$, respectively.

Two integer matrices ${\sf M}, {\sf N}$ are \emph{equivalent} if there exist unimodular matrices $\sf P$ and $\sf Q$ with entries in $\mathbb{Z}$ satisfying ${\sf M}={\sf PNQ}$.
In such case, this relation is denoted by ${\sf N}\sim{\sf M}$.
Therefore, if $\sf M$ and $\sf N$ are equivalent, then $\sf M$ can be transformed into $\sf N$ by means of the following operations:
\begin{enumerate}
  \item swap any two rows or any two columns.
  \item add an integer multiple of one row to another row.
  \item add an integer multiple of one column to another column.
  \item multiply any row or column by $\pm 1$.
\end{enumerate}
The \emph{Smith normal form} of a integer matrix $\sf M$, denoted by $\SNF(\sf M)$, is the unique diagonal matrix $\diag(f_1,\dots,f_r,0,\dots,0)$ equivalent to $\sf M$ such that $r=rank(\sf M)$ and $f_i|f_j$ for $i<j$.
The \emph{invariant factors} (or \emph{elementary divisors}) of $\sf M$ are the integers in the diagonal of the $\SNF(\sf M)$.
The reader might recall that in {\tt sagemath} the function {\tt elementary\_divisors()} returns the invariant factors of a integer matrix.

By considering an $m\times n$ matrix $\sf M$ with integer entries as a linear map ${\sf M}:\mathbb{Z}^n\rightarrow \mathbb{Z}^m$, the {\it cokernel} of $\sf M$ is the quotient module $\mathbb{Z}^{m}/{\rm Im}\, \sf M$.
It is known that if ${\sf N}\sim {\sf M}$, then $coker({\sf M})=\mathbb{Z}^n/{\rm Im} {\sf M}\cong\mathbb{Z}^n/{\rm Im} {\sf N}=coker(\sf N)$.
Therefore, as the fundamental theorem of finitely generated Abelian groups states, the cokernel of $\sf M$ can be described as:
$coker({\sf M})\cong \mathbb{Z}_{f_1} \oplus \mathbb{Z}_{f_2} \oplus \cdots \oplus\mathbb{Z}_{f_{r}} \oplus \mathbb{Z}^{n-r}$,
where $r=\rank({\sf M})$, and $f_1, f_2, \dots, f_{r}$ are the {\it invariant factors} of $\sf M$.
Aside to the homology groups, this finitely generated Abelian group becomes a graph invariant when we take the matrix $\sf M$ to be a matrix associated with the graph, say, the adjacency or Laplacian matrices.
Recall that the cokernel of the adjacency matrix is known as the {\it Smith group}, and the torsion part of the cokernel of the Laplacian matrix is known as the {\it critical group} or {\it sandpile group}.

Few is known on the SNF of distance matrices, some of the results obtained might be found in \cite{MR4884117,at,az,MR4110282,MR3615533,MR3524483,HW}.
Graham-Lovász-Pollak \cite{GL,GP} obtained the celebrated formula $$\det({\sf D}(T_{n+1}))=(-1)^nn2^{n-1},$$ for any tree $T_{n+1}$ with $n+1$ vertices.
Later, Hou and Woo \cite{HW} extended this formula to the SNF obtaining that $\SNF({\sf D}(T_{n+1}))={\sf I}_2\oplus 2{\sf I}_{n-2}\oplus [2n]$, for any tree $T_{n+1}$ with $n+1$ vertices.
There have been few more attempts in extending Graham-Lovász-Pollak formula, see \cite{alfaro2025characterizationgraphstrivialdistance,MR4818661,MR3891112,MR4623916,MR4679932,MR4020870,MR4209013,MR4103841,MR2231094} for instance.


The aim of this paper is to prove the following formula for the SNF of the $k$-distance matrix of any $k$-tree with $n\geq k+2$ vertices
$$\SNF({\sf D}^k(T_n))={\sf I}_{(k-1)(n-k)+2}\oplus (k+1){\sf I}_{n-k-2}\oplus [k(k+1)(n-k)].$$ 
As by product, we also obtain the determinant of ${\sf D}^k(T_n)$. 
These results are interesting since there are 2-trees with the same number of vertices whose 1-distance matrices have different SNF and determinant, this was already observed in \cite{MR3891112}.

\section{Results}


There is only one $k$-tree which consists of a $k$-clique, for this reason the $k$-distance matrix is the one-entry zero matrix.
And, there is only one $k$-tree $T_{k+1}$ consisting of a $(k+1)$-clique, therefore each of the $(k+1)$ $k$-cliques are adjacent to each other, thus ${\sf D}^k(T_{k+1})={\sf J}_{k+1}-{\sf I}_{k+1}$. 
It is not difficult to see that $\SNF({\sf J}_{k+1}-{\sf I}_{k+1})=\diag(1,\dots,1,k)$.
Therefore, we have the following result.

\begin{lemma}
    The Smith normal form of $k$-distance matrix of the unique $k$-tree with $k+1$ vertices is $\diag(1,\dots,1,k)$.
\end{lemma}



Let $T_{k+1}$ be the unique $k$-tree with $k+1$ vertices and its $k$-cliques labeled $1,\ldots,k+1$. 
For $i\in \{1,\ldots,k+1\}$, let $\alpha_{i}(T_{k+1})$ be the $k$-tree $T_{k+2}$ resulting from attaching a $(k+1)$-clique to $T_{k+1}$ through the $i$-{\it th} $k$-clique and labeling the $k$ new $k$-cliques from $k+2$ to $2k+1$. 
Inductively, every $k$-tree $T_n$ with $n$ vertices and its labeled $k$-cliques is the result of the application of $n-k-1$ successive analogous attachments, that is, $T_n=\alpha_{i_{n-k-1}}(\cdots\alpha_{i_1}(T_{k+1})\cdots)$, where $i_1,\dots,i_{n-k-1}$ are the indices of the $k$-cliques selected to attach the $(k+1)$-clique at each step. 
The {\it order} $s(n)$ of the matrix ${\sf D}^{k}(T_{n})$ is $k(n-k)+1$, since $k$ new $k$-cliques were added to $T_{k+1}$ with the application of each attachment $\alpha_{i_j}$. 
We assume the $k$-cliques of $T_n$ are labeled so that the entry ${\sf D}_{ij}^{k}(T_n)$ is the distance between the $i$-{\it th} $k$-clique and the $j$-{\it th} $k$-clique. 

Let $T_{n+1}=\alpha_i(T_n)$. 
Since the distance between the $k$-cliques of $T_n$ do not change after applying $\alpha_i$, then we assume ${\sf D}^{k}(T_n)$ is the upper left block of order $s(n)$ of ${\sf D}^{k}(T_{n+1})$. 
For $l\in\{1,\dots,k\}$, the distance from the $(s(n)+l)$-{\it th} $k$-clique to the $j$-{\it th} clique with $j\in\{1,\dots,s(n)\}$ is ${\sf D}_{ji}(T_n)+1$, since only one more $(k+1)$-clique was attached to $T_{n}$ through the $i$-{\it th} clique. 
Thus, if $d_{ij}={\sf D}^{k}_{ij}(T_{n})$, then ${\sf D}^{k}(T_{n+1})$ is equal to
$$\begin{bmatrix}
d_{11} & d_{12}  &\ldots & d_{1s(n)} & d_{1i}+1 & d_{1i}+1 &\cdots & d_{1i}+1\\
d_{21} & d_{22}  &\ldots & d_{2s(n)} & d_{2i}+1 &d_{2i}+1 &\cdots & d_{2i}+1\\
\vdots &\vdots &\cdots & \vdots& \vdots&\vdots&\cdots & \vdots\\
d_{s(n)1} & d_{s(n)2}  &\ldots & d_{s(n)s(n)} & d_{s(n)i}+1 &d_{s(n)i}+1 &\cdots & d_{s(n)i}+1\\
d_{i1}+1 & d_{i2}+1  &\ldots & d_{is(n)}+1 & 0 & 1 &\cdots & 1\\
d_{i1}+1 & d_{i2}+1  &\ldots & d_{is(n)}+1 & 1 & 0 &\cdots & 1\\
\vdots &\vdots &\cdots & \vdots& \vdots&\vdots&\cdots & \vdots\\
d_{i1}+1 & d_{i2}+1  &\ldots & d_{is(n)}+1 & 1 & 1 &\cdots & 0\\
\end{bmatrix}.
$$

Given ${\sf B}\in M_{r}(\mathbb{Z})$ and $i\in\{1,\dots, r\}$, let ${\sf A}_i({\sf B})\in M_{r+k}(\mathbb{Z})$ be the block matrix
$$
\begin{bmatrix}
{\sf B} & {\sf C}_{i} \\
{\sf C}_i^{T} & {\sf J}_k-{\sf I}_k \\
\end{bmatrix},
$$
where each column of ${\sf C}_i$ is the $i$-{\it th} column of $\sf B$ plus a column vector of ones. 
From the previous description of $\alpha_i$, if $T_{n+1}=\alpha_i(T_n)$, then
\begin{equation}
{\sf A}_i({\sf D}^{k}(T_{n}))={\sf D}^k(\alpha_i(T_{n})), \qquad \text{ for } i\in\{1,\dots,s(n)\}.
\end{equation}
In particular, if $n=k+1$, then for the unique $k$-tree with $k+1$ vertices, we have ${\sf A}_i({\sf D}^{k}(T_{k+1}))={\sf D}^k(\alpha_i(T_{k+1})),$ for any $i\in\{1,\dots,k+1\}$. 
We have observed every $k$-tree with $n$ vertices is the result of $n-k-1$ compositions: $T_n=\alpha_{i_{n-k-1}}(\cdots\alpha_{i_1}(T_{k+1})\cdots)$, thus, inductively,
\begin{equation}\label{eqindAi}
{\sf D}^{k}(T_{n})={\sf A}_{i_{n-k-1}}(\cdots {\sf A}_{i_1}({\sf D}^{k}(T_{k+1})\cdots).
\end{equation}
The way in which ${\sf D}^{k}(T_{n})$ is written depends on the labels of the $k$-cliques in $T_{n}$, however, we can always write ${\sf D}^{k}(T_n)$ as in Equation~\ref{eqindAi}.

\begin{lemma}\label{LemaEquivRho}
Let $T_{n}$ be a $k$-tree with $n$ vertices and its $k$-cliques labeled from 1 to $s(n)$. 
If $\rho$ is a permutation of the set of labels and $T^{\rho}_{n}$ is the tree $T_n$ with its $k$-cliques labeled $\rho(1),\ldots, \rho(s(n))$, then ${\sf D}^{k}(T^{\rho}_{n})$ is equivalent to ${\sf D}^{k}(T_{n})$.
\end{lemma}
\begin{proof} 
The $k$-cliques labeled $i$ and $j$ in $T^{\rho}_{n}$ are labeled $\rho^{-1}(i)$ and $\rho^{-1}(j)$ in $T_{n}$. 
So, 
\begin{equation}\label{DTDTr}
{\sf D}^{k}_{ij}(T^{\rho}_{n})={\sf D}^{k}_{\rho^{-1}(i)\rho^{-1}(j)}(T_{n}).
\end{equation}
Let $\sf P$ be the matrix whose $(i,j)$-entry is the Kronecker delta $\delta_{\rho^{-1}(i)j}$. 
Since $\sf P$ is a permutation matrix, then ${\sf P}^{-1}_{ij}=\delta_{i\rho^{-1}(j)}$ and $\det({\sf P})=\det({\sf P}^{-1})=\pm1$. 
Now
\begin{align*}
({\sf PD}^{k}(T_{n}){\sf P}^{-1})_{ij} &=\sum_{l=1}^{s(n)}({\sf PD}^{k}(T_{n}))_{il}\delta_{l\rho^{-1}(j)}=({\sf PD}^{k}(T_{n}))_{i\rho^{-1}(j)}\\
&=\sum_{l=1}^{s(n)}\delta_{\rho^{-1}(i)l}{\sf D}^{k}_{l\rho^{-1}(j)}(T_{n})={\sf D}^{k}_{\rho^{-1}(i)\rho^{-1}(j)}(T_{n}).
\end{align*}
The result turns out from Equation~\ref{DTDTr}.
\end{proof}

\begin{example}
In Figure~\ref{fig:labelingtrees}, there are two different labelings of the $2$-cliques of a $2$-tree. 
\begin{figure}[!ht]
    \centering
    \begin{tabular}{c@{\extracolsep{1cm}}c}
    
\begin{tikzpicture}
    \coordinate (A) at (0,0);
    \coordinate (B) at (2,0);
    \coordinate (C) at (1,1.5);
    \coordinate (D) at (1,-1.5);

    \draw[thick] 
        (A) -- (B) node[midway, above] {$4$}
        -- (C) node[midway, right] {$5$}
        -- (A) node[midway, left] {$1$};
    \draw[thick]
        (A) -- (D) node[midway, left] {$3$}
        -- (B) node[midway, right] {$2$};
\end{tikzpicture}
&
\begin{tikzpicture}
    \coordinate (A) at (0,0);
    \coordinate (B) at (2,0);
    \coordinate (C) at (1,1.5);
    \coordinate (D) at (1,-1.5); 

    \draw[thick] 
        (A) -- (B) node[midway, above] {$1$}
        -- (C) node[midway, right] {$3$}
        -- (A) node[midway, left] {$2$};

    \draw[thick]
        (A) -- (D) node[midway, left] {$4$}
        -- (B) node[midway, right] {$5$};
\end{tikzpicture}
\end{tabular}
\caption{Two labelings of the 2-cliques of a $2$-tree.}
\label{fig:labelingtrees}
\end{figure}
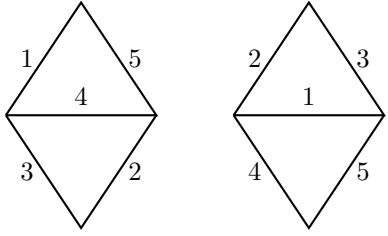
Let $T_5$ be the one on the right side and let $T_5^{\rho}$ be the one on the left side, with $\rho$ defined by the labels on the corresponding $2$-cliques. 
The following equation shows the equivalence of Lemma~\ref{LemaEquivRho} for ${\sf D}^2(T_5^{\rho})$ and ${\sf D}^2(T_5)$: 
$$\begin{bmatrix}
0 & 2 & 2 & 1 & 1\\
2 & 0 & 1& 1 & 2\\
2 & 1 & 0 & 1& 2\\
1 & 1 & 1& 0 & 1\\
1 & 2 & 2 & 1 & 0\\
\end{bmatrix} = 
\begin{bmatrix}
0 & 1 & 0 & 0 & 0\\
0 & 0 & 0& 0& 1\\
0 & 0 & 0 & 1& 0\\
1 & 0 & 0& 0 & 0\\
0 & 0 & 1 & 0 & 0\\
\end{bmatrix}
\begin{bmatrix}
0 & 1 & 1 & 1 & 1\\
1 & 0 & 1& 2 & 2\\
1 & 1 & 0 & 2& 2\\
1 & 2 & 2& 0 & 1\\
1 & 2 & 2 & 1 & 0\\
\end{bmatrix}
\begin{bmatrix}
0 & 0 & 0 & 1 & 0\\
1 & 0 & 0& 0 & 0\\
0 & 0 & 0 & 0& 1\\
0 & 0 & 1& 0 & 0\\
0 & 1 & 0 & 0 & 0\\
\end{bmatrix}.
$$
\end{example}

From now on, we assume ${\sf D}^{k}(T_n)$ is written as in Equation~$\ref{eqindAi}$ for any $k$-tree with $n$ vertices.

\begin{theorem}\label{TheoEquiv}
For $k$ and $n$ such that $k\geq 1$ and $n\geq k+2$, the $k$-distance matrices 
of all $k$-trees with $n$ vertices are equivalent.
\end{theorem}
\begin{proof}
Let $T_n$ be a $k$-tree with $n$ vertices.
From Equation~$\ref{eqindAi}$, we have ${\sf D}^{k}_{ij}(T_n)={\sf D}^{k}_{ji}(T_n)={\sf D}_{i_{n-k-1}j}+1$, for $i$ and $j$ such that $1\leq j\leq s(n-1)$ and $s(n-1)+1\leq i\leq s(n)$. 
In particular, ${\sf D}^{k}_{ii_{n-k-1}}(T_n)={\sf D}^{k}_{i_{n-k-1}i}(T_n)=1$, for $i$ such that $s(n-1)+1\leq i\leq s(n)$. 
Thus, subtracting the $i_{n-k-1}$-{\it th} column to the last $k$ columns of ${\sf D}^{k}(T_n)$ and the $i_{n-k-1}$-{\it th} row to the last $k$ rows of ${\sf D}^{k}(T_n)$, we obtain the following equivalent matrix, where $d_{ij}={\sf D}^{k}_{ij}(T_n)$,
$$\begin{bmatrix}
d_{11} & d_{12}  &\ldots & d_{1s(n-1)} & 1 & 1 &\cdots & 1\\
d_{21} & d_{22}  &\ldots & d_{2s(n-1)} & 1 &1 &\cdots & 1\\
\vdots &\vdots &\ddots & \vdots& \vdots&\vdots&\ddots & \vdots\\
d_{s(n-1)1} & d_{s(n-1)2}  &\ldots & d_{s(n-1)s(n-1)} & 1 &1 &\cdots & 1\\
1 & 1  &\ldots & 1 & -2 & -1 &\cdots & -1\\
1 & 1  &\ldots & 1 & -1 & -2 &\cdots & -1\\
\vdots &\vdots &\ddots & \vdots& \vdots&\vdots&\ddots & \vdots\\
1 & 1  &\ldots & 1 & -1 & -1 &\cdots & -2\\
\end{bmatrix}.
$$
Similarly, ${\sf D}^{k}_{ij}(T_n)={\sf D}^{k}_{ji}(T_n)={\sf D}_{i_{n-k-2}j}(T_n)$+1, for $i$ and $j$ such that $1\leq j\leq s(n-1)$ and $s(n-2)+1\leq i\leq s(n-1)$. 
In particular, ${\sf D}^{k}_{ii_{n-k-2}}(T_n)={\sf D}^{k}_{i_{n-k-2} i}(T_n)=1$, for $i$ such that $s(n-2)+1\leq i\leq s(n-1)$. 
Subtracting the $i_{n-k-2}$-{\it th} column to the columns $s(n-2)+1,\dots, s(n-1)$ of ${\sf D}^{k}(T_n)$ and the $i_{n-k-2}$-{\it th} row to the rows $s(n-2)+1,\dots, s(n-1)$ of ${\sf D}^{k}(T_n)$, we obtain the following equivalent matrix
$$\left[\begin{array}{cccccccccccc}
d_{11} & d_{12}  &\ldots & d_{1s(n-2)} & 1 & 1 & \cdots & 1 & 1 & 1 & \cdots& 1\\
d_{21} & d_{22}  &\ldots & d_{2s(n-2)} & 1 & 1 & \cdots & 1 & 1 & 1 &\cdots & 1\\
\vdots & \vdots  &\ddots & \vdots      & \vdots &\vdots&\ddots & \vdots&\vdots&\vdots&\ddots & \vdots\\
d_{s(n-2)1} & d_{s(n-2)2}  &\ldots & d_{s(n-2)s(n-2)} & 1 &1 &\cdots & 1 & 1 &1 &\cdots & 1\\
1 & 1  &\ldots & 1 & -2 &-1 &\cdots & -1 & 0 & 0 &\cdots & 0\\
1 & 1  &\ldots & 1 & -1 &-2 &\cdots & -1 & 0 & 0 &\cdots & 0\\
\vdots &\vdots &\ddots & \vdots& \vdots&\vdots&\ddots & \vdots& \vdots&\vdots&\ddots & \vdots\\
1 & 1  &\ldots & 1 & -1 &-1 &\cdots & -2 & 0 & 0 &\cdots & 0\\
1 & 1  &\ldots & 1 & 0 &0 &\cdots & 0 & -2 & -1 &\cdots & -1\\
1 & 1  &\ldots & 1 & 0 &0 &\cdots & 0 & -1 & -2 &\cdots & -1\\
\vdots &\vdots &\ddots & \vdots& \vdots&\vdots&\ddots & \vdots& \vdots&\vdots&\ddots & \vdots\\
1 & 1  &\ldots & 1 & 0 &0 &\cdots & 0 & -1 & -1 &\cdots & -2\\
\end{array}
\right].
$$
Continuing with the elementary operations on all the $k\times k$ blocks, we obtain
$$\footnotesize
\left[\begin{array}{cccccccccccccccccccc}
0 & 1  & \cd & 1 & 1 & 1 & \cd & 1 & 1 & 1 & \cd& 1& 1 & 1 & 1& 1 & 1 & \cd& 1\\
1 & 0  & \cd & 1 & 1 & 1 & \cd & 1 & 1 & 1 & \cd& 1& 1 & 1 & 1& 1 & 1 & \cd& 1\\
\vd & \vd  &\dd & \vd & \vd &\vd &\dd &\vd & \vd&\vd&\dd &\vd & \vd &\vd &\dd & \vd &\vd &\dd &\vd\\
1 & 1  & \cd & 0 & 1 & 1 & \cd & 1 & 1 & 1 & \cd & 1& 1 & 1 & 1& 1 & 1 & \cd & 1\\
1 & 1  & \cd & 1 & 0 & 1 & \cd & 1 & 1 & 1 & \cd & 1& 1 & 1 & 1& 1 & 1 & \cd & 1\\ 
1 & 1  & \cd & 1 & -2 & -1 & \cd & -1 &  &  &  &  &  &  &  &  &  &  & \\
1 & 1  & \cd & 1 & -1 & -2 & \cd & -1 &  &  &  &  &  &  &  &  &  &  & \\
\vd & \vd  &\dd & \vd & \vd &\vd &\dd & \vd &  &  &  &  &  &  &  &  &  &  &  \\
1 & 1  & \cd & 1 & -1 & -1 & \cd & -2 &  &  &  &  &  &  &  &  &  &  & \\
1 & 1  & \cd & 1 &  &  &  &  &  -2 & -1 & \cd & -1 &  &  &  &  &  &  & \\
1 & 1  & \cd & 1 &  &  &  &  &  -1 & -2 & \cd & -1 &  &  &  &  &  &  & \\
\vd & \vd  &\dd & \vd &  &  &  &  & \vd & \vd & \dd & \vd  &  &  &  &  &  &  &  \\
1 & 1  & \cd & 1 &  &  &  &  & -1& -1 & \cd & -2 &  &  &  &  &  &  & \\
1 & 1  & \cd & 1 &  &  &  &  &   &  &  &  &  &   &  &  &  &  & \\
1 & 1  & \cd & 1 &  &  &  &  &  &  &  &  &  &  \ddots &  &  &  &  & \\
\vd & \vd  &\dd & \vd &  &  &  &  &  &  &  &  &   &  &  &  &  &  &  &  \\
1 & 1  & \cd & 1 &  &  &  &  &   &  &  &  &  &  &  & -2 & -1  & \cd & -1\\
1 & 1  & \cd & 1 &  &  &  &  &  &  &  &  &  &  &  & -1 &-2  & \cd & -1\\
\vd & \vd  &\dd & \vd &  &  &  &  &  &  &  &   &  &  &  & \vd & \vd & \dd & \vd \\
1 & 1  & \cd & 1 &  &  &  &  & &  &  &  &  &  &  & -1 &-1  & \cd & -2\\
\end{array}
\right].
$$

Finally, subtracting the first row to rows $2,\ldots,k+1$ and the first column to columns $2,\ldots,k+1$ we conclude ${\sf D}^{k}(T_n)$ is equivalent to
$$\footnotesize
\left[\begin{array}{cccccccccccccccccccc}
0 & 1  & 1 & \cd & 1 & 1 & 1 & \cd & 1 & 1 & 1& \cd& 1 & \cd &  & 1 & 1 & \cd& 1\\
1 & -2  & -1 & \cd & -1 &  &  &  &  &  & & &  &  & &  &  &  & \\
1 & -1  & -2 & \cd & -1 &  &  &  &  &  & & &  &  & &  &  &  & \\
\vd & \vd  &\vd & \dd & \vd &  &  &  &  &  &  &  &  &  &  &  &  &  &  \\
1 & -1  & -1 & \cd & -2  &  &  &  &  &  &  & &  &  &  &  &  &  & \\
1 &  &  &  &  & -2 & -1 & \cd & -1 &  &  &  &  &  &  &  &  &  &  \\
1 &  & &  &  & -1 & -2 & \cd & -1 &  &  &  &  &  &  &  &  &  &  \\
\vd &    &  &  & &\vd & \vd & \dd & \vd &  &  &  &  &  &  &  &  &  &\\
1 &  & &  &  & -1 & -1 & \cd & -2 &  &  &  &  &  &  &  &  &  &  \\
1 &  &  &  &  &  &  &  &  &  -2 & -1 & \cd & -1 &  &  &  &  &  &  \\
1 &  &  &  &  &  &  &  &  &  -1 & -2 & \cd & -1 &  &  &  &  &  &  &  \\
\vd  & &  &   & &  &  &  &  & \vd & \vd & \dd & \vd  &  &  &  &  &  &  & \\
 & &  &  & &  &  &  &  & -1& -1 & \cd & -2 &  &  &  &  &  &  \\
\vd &   &  &  &  &  &  &  &   &  &  &  &  &   &  &  &  &  & \\
 &  &  & &  &  &  &  &  &  &  &  &  &  \ddots &  &  &  &  & \\
\vd &   & &  &  &  &  &  &  &  &  &  &   &  &  &  &  &  &  &  \\
1 &   &  &  &  &  &  &  &   &  &  &  &  &  &  & -2 & -1  & \cd & -1\\
1 &   &  &  &  &  &  &  &  &  &  &  &  &  &  & -1 &-2  & \cd & -1\\
\vd &   & &  &  &  &  &  &  &  &  &   &  &  &  & \vd & \vd & \dd & \vd \\
1 &   &  &  &  &  &  &  & &  &  &  &  &  &  & -1 &-1  & \cd & -2\\
\end{array}
\right]
$$

\end{proof}

\begin{lemma}\label{LemaMatrixM}
For $k\geq 2$ and ${\sf M}_k=-{\sf J}_{k}-{\sf I}_{k}$,   
$\SNF({\sf M}_k)=  \diag(1,\ldots,1,k+1)$.
\end{lemma}

\begin{proof}
Consider the matrices
$${\sf P_M}=\begin{bmatrix}
1 & 0  &\ldots & 0 & 0\\
0 & 1  & \ldots & 0 & 0\\
\vdots &\vdots &\ddots & \vdots& \vdots\\
0 & 0 &\ldots & 1 & 0\\
-k & -k &\ldots & -k & 1\\
\end{bmatrix}
\text{ and }
{\sf Q_M}=\begin{bmatrix}
0 & 1 & \ldots & 1 & 1\\
1 & 0 &  \ldots & 1 & 1\\
\vdots &\vdots &\ddots & \vdots& \vdots\\
1 & 1 &\ldots & 0 & 1\\
-(k-1) & -(k-1) &\ldots & -(k-1) & -k\\
\end{bmatrix}.$$
Using elementary row operations, we obtain ${\sf P_M}\sim {\sf I}_{k}$  and ${\sf Q_M}\sim -{\sf I}_k$. 
Thus 
\begin{equation}\label{detsPMQM}
\det({\sf P_M})=1\qquad\text{ and }\qquad \det({\sf Q_M})=(-1)^k.
\end{equation}
A direct computation shows ${\sf P_M M}_k {\sf Q_M}=\diag(1,\dots,1,k+1)$.
\end{proof}

Let us recall the following result of \cite{HW}.

\begin{lemma}\cite[Corollary 2]{HW}\label{lema:HouWoo}
    The Smith normal form of the matrix
    \[
    \begin{bmatrix}
        c & b{\sf 1}^T\\
        b{\sf 1} & \diag(a,\dots,a)
    \end{bmatrix}
    \]
    of order $n+1$ is $\diag(\Delta_1,\Delta_2,a,\dots,a,\Delta_{n+1})$,
    where
    \[
    \Delta_1=\gcd(a,b,c), \Delta_2=\frac{\gcd(a^2,b^2,ca,,ba)}{\gcd(a,b,c)} \text{ and } \Delta_{n+1}=\frac{ac-nb^2}{\gcd(a,b,c)}.
    \]
\end{lemma}

\begin{theorem}
Let $k\geq 1$ and $n\geq k+2$.
For any $k$-tree $T_n$ with $n$ vertices,
$$\SNF({\sf D}^{k}(T_n))={\sf I}_{(k-1)(n-k)+2}\oplus (k+1){\sf I}_{n-k-2}\oplus [k(k+1)(n-k)].$$
\end{theorem}
\begin{proof}
According to Theorem~\ref{TheoEquiv}, we only need to compute the Smith normal form of the matrix
$$\begin{bmatrix}
0 & {\bf 1}^T & {\bf 1}^T & \cd & {\bf 1}^T \\
{\bf 1} & {\sf M}_k &  &  &  \\
{\bf 1} &  & {\sf M}_k &  &  \\
\vd &  &  & \ddots & \\
{\bf 1} &  &  &  & {\sf M}_k \\
\end{bmatrix},$$
where ${\sf M}_k$ is the matrix defined in Lemma~\ref{LemaMatrixM}, which appears $n-k$ times. 
Note that in the special case $k=1$, instead of having the block matrices ${\sf M}_k$, the matrix has an entry $-2$.

Using Lemma~\ref{LemaMatrixM}, the product
$$
\begin{bmatrix}
1 & {\bf 0}^T & {\bf 0}^T & \cd & {\bf 0}^T \\
{\bf 0} & {\sf P_M} &  &  &  \\
{\bf 0} &  & {\sf P_M} &  &  \\
\vd &  &  & \ddots & \\
{\bf 0} &  &  &  & {\sf P_M} \\
\end{bmatrix}
\begin{bmatrix}
0 & {\bf 1}^T & {\bf 1}^T & \cd & {\bf 1}^T \\
{\bf 1} & {\sf M}_k &  &  &  \\
{\bf 1} &  & {\sf M}_k &  &  \\
\vd &  &  & \ddots & \\
{\bf 1} &  &  &  & {\sf M}_k \\
\end{bmatrix}
\begin{bmatrix}
1 & {\bf 0}^T & {\bf 0}^T & \cd & {\bf 0}^T \\
{\bf 0} & {\sf Q_M} &  &  &  \\
{\bf 0} &  & {\sf Q_M} &  &  \\
\vd &  &  & \ddots & \\
{\bf 0} &  &  &  & {\sf Q_M} \\
\end{bmatrix}$$
is the matrix
$$\begin{bmatrix}
0 & -{\bf 1}^T & -{\bf 1}^T & \cd & -{\bf 1}^T \\
{\bf v} & {\sf D_M} &  &  &  \\
{\bf v} &  & {\sf D_M} &  &  \\
\vd &  &  & \ddots & \\
{\bf v} &  &  &  & {\sf D_M} \\
\end{bmatrix},$$
where ${\sf D_M}=\diag(1,\dots, 1, k+1)$ and 
$\bf v$ is the column vector whose entries are the sums of the rows of $\sf P_M$, that is, ${\bf v}^T=(1,\dots,1,-k(k-1)+1)$.

If ${\bf w}^{T}=(-2,\ldots,-2,k-2)$, then
$$
\begin{bmatrix}
0 & -{\bf 1}^T & -{\bf 1}^T & \cd & -{\bf 1}^T \\
{\bf v} & {\sf D_M} &  &  &  \\
{\bf v} &  & {\sf D_M} &  &  \\
\vd &  &  & \ddots & \\
{\bf v} &  &  &  & {\sf D_M} \\
\end{bmatrix}
\begin{bmatrix}
1 & {\bf 0}^T & {\bf 0}^T & \cd & {\bf 0}^T \\
{\bf w} & {\sf I} &  &  &  \\
{\bf w} &  & {\sf I} &  &  \\
\vd &  &  & \ddots & \\
{\bf w} &  &  &  & {\sf I} \\
\end{bmatrix}
=
\begin{bmatrix}
k(n-k) & -{\bf 1}^T & -{\bf 1}^T & \cd & -{\bf 1}^T \\
-{\bf 1} & {\sf D_M} &  &  &  \\
-{\bf 1} &  & {\sf D_M} &  &  \\
\vd &  &  & \ddots & \\
-{\bf 1} &  &  &  & {\sf D_M} \\
\end{bmatrix}.$$
Now, we can transform this last matrix to
$${{\sf I}_{(k-1)(n-k)}}\oplus\begin{bmatrix}
n-k & 1 & 1 & \cd & 1 \\
1 & k+1 &  &  &  \\
1 &  & k+1 &  &  \\
\vd &  &  & \ddots & \\
1 &  &  &  & k+1 \\
\end{bmatrix}.$$
This matrix has the form described in Lemma~\ref{lema:HouWoo} with $a=k+1$, $b=1$, $c=n-k$ and the order of the matrix is equal to $n-k+1$. 
Therefore, 
$$\SNF
\begin{bmatrix}
    n-k & {\sf 1}^T\\
    {\sf 1} & \diag(k+1,\dots,k+1)
\end{bmatrix}
=\diag(1,1,k+1,\ldots,k+1, k(k+1)(n-k)).$$
Thus the result follows.
\end{proof}

\begin{corollary}
$\det({\sf D}^{k}(T_n))=(-1)^{k(n-k)}k(k+1)^{n-k-1}(n-k)$.
\end{corollary}
\begin{proof}
The only term we need to justify is $(-1)^{k(n-k)}$. All the elementary operations described to obtain $\SNF({\sf D}^{k}(T_n))$ were performed with matrices of determinant one, except the multiplication by ${\sf Q_M}$. 
This matrix was employed $n-k$ times and from Equation~\ref{detsPMQM} has determinant $(-1)^k$.

In the special case $k=1$, instead of the matrix ${\sf M_k}$ we had a $-2$ term $n-k$ times. 
So in this case we have a factor of $(-1)^{n-k}$ instead, which corresponds with the factor of the determinant formula.
\end{proof}

\section*{Acknowledgement}
The research of C. A. Alfaro is partially supported by Sistema Nacional de Investigadoras e Investigadores grant number 220797.
The research of Jesús Uriel Medrano is partially supported by Secretaría de Ciencia, Humanidades, Tecnología e Innovación grant number 1343582.
The research of Iván Téllez Téllez is partially supported by Sistema Nacional de Investigadoras e Investigadores grant number 362694.

\bibliographystyle{plain}
\bibliography{bibliography}

\appendix

\section{Code to generate all $k$-trees with $n$ vertices}\label{sec:generation of k trees}
\begin{lstlisting}
def k_trees(k,nmax):
    G = [[graphs.CompleteGraph(k)]]
    for n in range(nmax-k):
        L = []
        for g in G[n]:
            for clique in sage.graphs.cliquer.all_cliques(g,k,k):
                h = g.copy()
                h.add_vertex()
                h.add_edges([(g.order(),v) for v in clique])
                bandera = True
                for l in L:
                    if h.is_isomorphic(l):
                        bandera = False
                        break
                if bandera:
                    L.append(h)
        G.append(L)
    return G
\end{lstlisting}

\section{Code to compute the $d$-distance matrix of a $k$-tree and its Smith normal form}\label{sec:k distance matrix}

\begin{lstlisting}
def d_distance_matrix(g,d):
    V = sorted(list(sage.graphs.cliquer.all_cliques(g,d,d)))
    print("\n" + str(d) + "-cliques")
    print(V)
    E = sorted(list(sage.graphs.cliquer.all_cliques(g,d+1,d+1)))
    print(str(d+1) + "-cliques")
    print(E)
    nV = len(V)
    nE = len(E)
    # Ak is the adjacency matrix of dimension k
    A = Matrix(nV,nV)
    for i in range(nV):
        for j in range(i+1,nV):
            flag = False
            for e in E:
                if set(V[i]).issubset(e) and set(V[j]).issubset(e):
                    flag = True
                    break
            if flag:
                A[i,j] = 1
                A[j,i] = 1
    H = Graph(A, format="adjacency_matrix")
    if H.is_connected():
        print("Distance matrix of dimension " + str(d))
        D = H.distance_matrix()
        print(D)
        print("SNF of the " + str(d) + "-dimensional distance matrix")
        print(D.elementary_divisors())
\end{lstlisting}

\end{document}